\documentclass[a4paper,11pt]{article}
\title{ Property (T) and Poincaré duality in dimension three}
\usepackage{geometry}

\date{}
 
\usepackage[T1]{fontenc}
\usepackage{comment}
\usepackage{microtype}
\usepackage{pinlabel}
\usepackage{float}
\usepackage[hypertexnames = false]{hyperref}
\usepackage{subfig}
\usepackage{bbm}
\usepackage{amssymb, amsthm,thmtools,cite,mathtools,setspace,enumitem,tikz-cd}
\usepackage{thm-restate}
\usepackage{bm}
\usepackage[T1]{fontenc} 
\usepackage{tikz-cd}

\usepackage{AlegreyaSans}
\usepackage[T1]{fontenc}
\usepackage[utopia,vvarbb]{newtxmath}
\usepackage{mathalfa}
\usepackage{cochineal}

\makeatletter
\newcommand\RedeclareMathOperator{%
	\@ifstar{\def\rmo@s{m}\rmo@redeclare}{\def\rmo@s{o}\rmo@redeclare}%
}
\newcommand\rmo@redeclare[2]{%
	\begingroup \escapechar\m@ne\xdef\@gtempa{{\string#1}}\endgroup
	\expandafter\@ifundefined\@gtempa
	{\@latex@error{\noexpand#1undefined}\@ehc}%
	\relax
	\expandafter\rmo@declmathop\rmo@s{#1}{#2}}
\newcommand\rmo@declmathop[3]{%
	\DeclareRobustCommand{#2}{\qopname\newmcodes@#1{#3}}%
}

\@onlypreamble\RedeclareMathOperator
\makeatother

\DeclareSymbolFont{bbold}{U}{bbold}{m}{n}

\newtheorem{thm}{Theorem}
\newtheorem{introthm}{Theorem}

\numberwithin{thm}{section}

\newtheorem{lem}[thm]{Lemma}

\newtheorem{prop}[thm]{Proposition}

\newtheorem{cor}{Corollary}

\DeclareMathOperator{\Fill}{\mathtt{fill}}

\DeclareMathOperator{\len}{{len}}

\newcommand{\Q}{{\mathbb Q}}
\newcommand{\Z}{{\mathbb Z}}

\newcommand{\R}{{\mathbb R}}

\newcommand{\cB}{\mathcal{B}}

\newcommand{\cL}{\mathcal{L}}

\newcommand{\cT}{\mathcal{T}}

\newcommand{\e}{\varepsilon}

\newcommand{\norm}[1]{\left\lVert#1\right\rVert}

\newcommand{\normalsub}{\unlhd}

\renewcommand{\d}{ \partial}

\definecolor{darkmidnightblue}{rgb}{0.0, 0.2, 0.4}
\hypersetup{
    colorlinks=true,
    linkcolor=darkmidnightblue,
    citecolor=darkmidnightblue,
    filecolor=darkmidnightblue,
    urlcolor=darkmidnightblue
}

\author{Cameron Gates Rudd}

\newcommand{\FP}{{\mathsf{FP}}}
\newcommand{\FL}{{\mathsf{FL}}}

\newcommand{\PD}{\mathsf{PD}}
\newcommand{\SL}{\mathsf{SL}}

\DeclarePairedDelimiter\floor{\lfloor}{\rfloor}

\begin{document}

\maketitle
\begin{abstract}
We use a recent result of Bader and Sauer on coboundary expansion to prove residually finite 3-dimensional Poincaré duality groups never have property (T). This implies such groups are never Kähler.
The argument applies to fundamental groups of (possibly non-aspherical) compact 3-manifolds as well, giving a new proof of a theorem of Fujiwara that states if the fundamental group of a compact 3-manifold has property (T), then that group is finite. The only consequence of geometrization needed in the proof is that 3-manifold groups are residually finite.
\end{abstract}


A countable discrete group $G$ has property (T) if $H^1(G;\pi)=0$ for all unitary representations $\pi$. We prove that this is incompatible with 3-dimensional Poincaré duality for residually finite groups.\footnote{
Note there do exist residually finite 3-dimensional Bieri-Eckmann duality groups with property (T), for instance, torsion free finite index subgroups of $\SL_3\Z$. Additionally, for $n>3$ it is known that there are non residually finite $\PD_n$ groups \cite{Mess}.
}

\begin{introthm}\label{introthm:PD3noT}
    Let $G$ be a residually finite $\PD_3$ group. Then $G$ does not have property (T).
\end{introthm}

The proof of Theorem \ref{introthm:PD3noT} uses recent work of Bader and Sauer
in which they discovered a new coboundary expansion phenomenon of groups with property (T) \cite{BaderSauer}.
 Theorem \ref{introthm:PD3noT} serves as proof of concept for a new geometric strategy for showing a group does not have property (T) using this higher coboundary expansion.

Theorem \ref{introthm:PD3noT} has the following corollary:

\begin{cor}\label{PD3noK}
    A residually finite $\PD_3$ group is not Kähler.
\end{cor}

This follows from the fact any $\PD_3$ Kähler group must have property (T), which is due to unpublished work of Delzant (see Theorem 4.3 in \cite{BiswasMjSeshadri}) or alternatively a theorem of Kotschick combined with work of Reznikov \cite{Reznikov,Kotschick}. 
Goldman and Donaldson, and independently Reznikov, conjectured that no 3-manifold group was Kähler \cite{Reznikov}. This conjecture was solved by Dimca and Suciu, and alternative proofs were given by Biswas-Mj-Seshadri and Kotschick \cite{DimcaSuciu,BiswasMjSeshadri,Kotschick}. Kotschick in fact proved rational $\PD_3$ groups with positive first Betti number cannot be Kähler, so this corollary extends Kotschick's theorem to residually finite (integral) $\PD_3$ groups with trivial first Betti number.

The argument used to prove Theorem \ref{introthm:PD3noT} also applies to the fundamental group of a compact 3-manifold (we do not assume the manifold is aspherical), recovering a theorem of Fujiwara \cite{Fujiwara}. Fujiwara's theorem is stated for geometric 3-manifolds and predates Perelman's proof of the geometrization conjecture, but applies to all compact 3-manifolds when combined with it \cite{PerelmanEntropy}.

\begin{introthm}[Fujiwara]\label{introthm:Manifold}
    Let $M$ be a compact 3-manifold whose fundamental group $G$ has property (T). Then $G$ is finite.
\end{introthm} 

We note that Wall conjectured that Poincaré duality groups of dimension 3 are exactly 3-manifold groups, however very little is known directly about 3-dimensional Poincaré duality groups. A related conjecture of Cannon says word hyperbolic groups with 2-sphere boundary contain finite index subgroups that are uniform lattices in $\text{SO}^+(3,1)$. Work of Bestvina and Mess shows that torsion free word hyperbolic groups with 2-sphere boundary are $\PD_3$ groups \cite{BestvinaMess}. Thus with this noted, Theorem \ref{introthm:PD3noT} and Corollary \ref{PD3noK} can be viewed as providing some evidence for these conjectures (with the added condition that the groups be residually finite). Note that it has been known since \cite{Wang} that lattices in $\text{SO}^+(n,1)$ do not have property (T).

Fujiwara's proof of Theorem \ref{introthm:Manifold} is already quite short, but requires the full strength of geometrization and appeals to other sources to rule out property (T) in the various cases that arise. The proof here uses much less geometrization input (only that 3-manifold groups are residually finite) and instead takes advantage of new results on property (T).

The key input is the following slight generalization of a theorem of Bader and Sauer. 

\begin{introthm}[$\FP_2$ version of Theorem 2.13 \cite{BaderSauer}]\label{introthm:FP2BaderSauerRexpansion} Let $G$ be a group of type $\FP_2$ with property (T). Let $F_*\to \Z$ be a length two partial resolution of finite rank free $\Z G$-modules with fixed bases. Endow $\hom_{\Z G}(F_i,\Z)$ with the $\ell^1$-norm $\norm{\cdot}$ induced by the bases. There is a constant $C$ depending on the partial resolution and bases such that for any finite index normal subgroup $H\normalsub G$ and every coboundary $\eta\in \hom_{\Z H}(F_2,\Z)$, there is a cochain $\omega\in  \hom_{\Z H}(F_1,\Z)$ such that $$d\omega = \eta \text{ and } \norm{\omega}\leq C \norm{\eta}.$$
\end{introthm}

Combined with Poincaré duality, the above theorem can be converted into a kind of codimension two linear isoperimetric inequality; see Proposition \ref{pdnhomologicalexpansion} for an exact statement in the present setting. Theorem 3.3 in \cite{BaderSauer} gives a geometric version for manifolds.
We note that in \cite{KielakKropholler}, Kielak and Kropholler showed that $n$-dimensional oriented Poincaré duality groups are either amenable or satisfy a linear homological isoperimetric inequality in codimension one. Proposition \ref{pdnhomologicalexpansion} combined with Section \ref{sec:hypandexpand} give a codimension two analogue for residually finite $\PD_n$ groups (of suitable finiteness type) with Property (T).
The following 2-dimensional nonexpansion result is then used to obtain a contradiction and prove Theorem \ref{introthm:PD3noT}.

\begin{introthm}\label{introthm:nonexpansion}
Let $M$ be a finite cell complex with residually finite infinite fundamental group $G.$ Let $G_i\normalsub G$ be a residual chain of subgroups. Set $M_i = \widetilde M/G_i.$ Then either $b_1(M_i)>0$ for some $i$, or for any constant $\e>0$,  there is a nontrivial boundary $z\in \d C_2(M_i;\Z)$ for some $M_i$, such that for any $2$-chain $A\in C_2(M_i;\Z)$ with $\d A = z$, one has $||z||< \e||A||.$
\end{introthm}

Note that the above result shows that a version of Bader and Sauer's theorem cannot hold for all groups that satisfy property $(\tau)$ with respect to some residual chain of subgroups. In particular, Brock-Dunfield and Boston-Ellenberg gave examples of closed hyperbolic 3-manifold groups that have property $(\tau)$ with respect to certain residual chains of subgroups, all of which have trivial first Betti number \cite{CDQHS,BE}. The above nonexpansion result combined with Poincaré duality shows these examples do not satisfy the uniform coboundary expansion of Theorem \ref{introthm:FP2BaderSauerRexpansion}. 

An interesting question is whether a stronger form of property $(\tau)$, which is still weaker than (T), implies coboundary expansion. With the (solved) 3-dimensional Lubotzky-Sarnak conjecture in mind, it is interesting to note that the above discussion shows $\PD_3$ groups could not satisfy this stronger version of ($\tau$).


\section{Coboundary expansion and property (T) }
In this section we outline Bader and Sauer's coboundary expansion result in the setting of groups of type $\FP_2$. The proofs follow from Bader and Sauer's work in \cite{BaderSauer}, with a slight argument needed to generalize their result on $\Z$-coefficients. The primary difference is that their work is in the context of cellular cochain complexes and here we work algebraically; note that in \cite{BaderSauer_survey}, a more algebraic account of related results is given. The purpose of this section is to set notation and convince the reader their proofs indeed generalize.

\subsection{Group cohomology and expansion}

Let $G$ be a group. Throughout, we work with left-modules over the group ring $\Z G$. We can convert left-modules into a right-modules and vice versa by twisting the action by the involution $g\mapsto g^{-1}$ of $G.$ If $M$ and $N$ are both left-modules, when we take their tensor product $M\otimes_{\Z G} N$, we use this twisted right-module structure on $M$. 

A partial $\Z G$-resolution of length $n$ consists of the first $n+1$ terms of a $\Z G$-resolution $F_*\to\Z$, where $\Z$ is the trivial module. A partial resolution is of type $\FP_n$ if it has length $n$ and every module $F_i$ is a finitely generated projective module. If instead every module $F_i$ is a finite rank free module, we say the partial resolution is of type $\FL_n$. 
A based partial $\Z G$-resolution of type $\FL_n$ additionally has a fixed $\Z G$-basis for each free module $F_i$ in the partial resolution.
If $G$ has a partial resolution of type $\FP_n,$ then $G$ is said to be of type $\FP_n.$

 We remark that every group of type $\FP_n$ has a partial resolution of type $\FL_n$; see \cite{Brown} Proposition VIII.4.3. We primarily work throughout with these based partial resolutions of type $\FL_n$ and refer to the groups as having type $\FP_n$.

The following basic fact will be useful.
\begin{lem}\label{overflow}
    Let $G$ be a group of type $\FP_2$. Let $g_1,\dots,g_m$ generate $G$ and let $X^1$ be the Cayley graph of this generating set. Then there is a finite (but typically incomplete) collection of relations such that after attaching 2-cells along the free $G$-orbits of these relations, the resulting 2-complex $X$ has trivial first homology and $G$ acts freely, cellularly, and cocompactly on $X.$ The augmented chain complex $C_*(X;\Z)$ with the $\Z G$-module structure induced by the $G$-action on $X$ and the basis given by choosing one cell for every free $G$-orbit is a based partial $\Z G$-resolution of type $\FL_2$, which we call the cellular partial resolution.
    \end{lem}

 Fix a based partial $\Z G$-resolution $F_*\to\Z$. Let $V$ be a normed $\Z G$-module. We define the chain complex $$C_*(F_*;V) := F_*\otimes_{\Z G}V$$ and cochain complex $$C^*(F_*;V):=\hom_{\Z G}(F_*,V).$$

If the coefficient module $V$ has the structure of a Banach space and the action of $G$ on $V$ is continuous, then the cochain complex $$\cdots \to C^i(F_*;V)\to C^{i+1}(F_*;V)\to\cdots$$ is a sequence of Fréchet spaces with continuous coboundary maps; the topology these complexes induce on the cohomology $H^i(G;V)$ is independent of the resolution.

When $F_i$ is a based finite rank free $\Z G$-module, there is an identification $C^i(F_*;V)\cong V^{k_i}$ and $C^i(F_*;V)$ inherits a norm. For $\bar v = (v_j)_1^{k_i}\in V^{k_i},$ the induced norm is given by $$||\bar v|| = \sum_{j=1}^{k_i} \norm{v_j},$$ where the norm $||v_j||$ is the norm on $V.$

If $H\normalsub G$, then because $\Z G$ is a free $\Z H$-module, the partial resolution $F_*\to \Z$ also gives a partial $\Z H$-resolution over the trivial $\Z H$-module $\Z$. If $F_i$ is finite rank as a $\Z G$-module and $H$ has finite index in $G,$ then $F_i$ is also a finite rank free $\Z H$-module.
We define $C_i(F_*|_H;V)$ and $C^i(F_*|_H;V)$ as above using the partial resolution $F_*\to\Z$ with the restriction $\Z H$-module structure.

For $H\normalsub G$ of finite index $m=[G:H]$, there is a (non-canonical) decomposition $$\Z G = \bigoplus\limits_{Hg \in H\backslash G} \Z H$$ inducing a decomposition $$C^i(F_*|_H;V) \cong \bigoplus_{Hg\in H\backslash G} V^{k_i}\cong V^{mk_i}.$$ The above decomposition defines a basis for $C^i(F_*|_H;V)$ associated to the initial free $\Z G$-module bases and choice of coset representatives; the norm induced by this basis is independent of the choices involved when $G$ acts on $V$ by linear isometries.

\begin{prop} Let $G$ be a group and $H\normalsub G$ a finite index subgroup. Let $V$ be a normed $\Z G$-module with isometric linear $G$ action. Let $F_*\to\Z$ be a based partial resolution of type $\FL_2$ and let $C^i(F_*|_H;V)$ be the associated cochain complex. The $\ell^1$-norm on $C^i(F_*|_H;V)$ induced by a basis coming from choice of coset representatives is independent of this choice. 
\end{prop}
\begin{proof}
    The norm on  $C^i(F_*|_H;V)$ is determined by the norm on $V$ and the basis, which depends only on the choices of coset representatives and the initial basis for $F_*$ as a $\Z G$-module. Choosing different coset representatives changes the decomposition of $C^i(F_*|_H;V)$ by a composition of factor permutations and the group action, which by hypothesis is isometric. As these are isometries of the $\ell^1$-norm, the $\ell^1$-norm is independent of the coset representative choice.
\end{proof}
Throughout this note, we work with normed $\Z G$-modules with isometric linear actions. After fixing a basis, we assume the norms we use are those described here.

A special class of normed coefficient modules is given by abstract $L$ spaces; a notion that generalizes spaces like $L^1(\Omega)$ for a measure space $\Omega$. Of relevance here is a special subclass denoted $\cL$, the specifics of this subclass are unimportant for the statements here, so we refer the reader to Section 1.1 of \cite{BaderSauer} for a discussion tailored to the present application. 
The fundamental result of Bader and Sauer is the following theorem about the induced topology on cohomology with coefficients in an abstract $L$ space in the class $\cL$ with an isometric and linear action.

\begin{thm}[$\FP_2$ version of Theorem 1.6 in \cite{BaderSauer}]
    Let $G$ be a group of type $\FP_2$ with property (T). Let $F_*\to \Z$ be a based partial $\Z G$-resolution of type $\FL_2$. Then for any abstract $L$-space $V$ in the class $\cL$ with isometric linear $G$ action, the cohomology $H^2(G,V)$ is Hausdorff.
\end{thm}
\begin{proof}
This follows from Corollary D in \cite{BaderGelanderMonod} combined with Lemma 29 in \cite{BaderSauer_survey} and the fact the class $\cL$ is closed under ultrapowers. \end{proof}

Using this, Bader and Sauer derive two expansion results.
The argument in \cite{BaderSauer} obtains these as a consequence of the Hausdorffness of $H^2(G,V)$ and properties of the class $\cL$ and its relation to property (T). In particular, it does not make any use of the cell structure, so applies verbatim in this setting as well.

\begin{thm}[$\FP_2$ version of Theorem 1.7 in \cite{BaderSauer}]
 Let $G$ be a group of type $\FP_2$ with property (T). Let $F_*\to \Z$ be a based partial $\Z G$-resolution of type $\FL_2$. Let $V$ be an abstract $L$-space in the class $\cL$ with an isometric linear $G$-action and endow the cochain complex with the $\ell^1$-norm induced by the bases and the norm on $V$. There is a constant $C$ depending on the partial resolution and bases such that for any coboundary $\eta\in C^2(F_*;V),$ there is a cochain $\omega\in C^1(F_*;V)$ such that $$d\omega=\eta \text{ and } ||\omega||\leq C||\eta||.$$
 \end{thm}

The main application of this result is the following.
\begin{thm}[$\FP_2$ version of Theorem 1.8 \cite{BaderSauer}]\label{FP2BaderSauerRexpansion} Let $G$ be a group of type $\FP_2$ with property (T). Let $F_*\to \Z$ be a based partial $\Z G$-resolution of type $\FL_2$. There is a constant $C$ depending on the partial resolution and bases such that for any finite index normal subgroup $H\normalsub G$ and every coboundary $\eta\in C^2(F_*|_H,\R)$,  there is a cochain $\omega\in C^1(F_*|_H,\R)$ such that $$d\omega = \eta \text{ and } \norm{\omega}\leq C \norm{\eta}.$$
\end{thm}

We refer to the constant $C$ above as the expansion constant of the cochain complex; these theorems give uniform bounds on expansion constants.

Bader and Sauer also prove that in the case the resolution comes from a cellular classifying space, one can upgrade the result above from real coefficients to integral coefficients. Their argument applies to the cellular free resolution of Lemma \ref{overflow} as their result is just about coboundary maps of 2-dimensional cell complexes.

\begin{thm}[Theorem 2.12 \cite{BaderSauer}]\label{linearprogram} Let $M$ be a finite cell complex such that $H^1(M;\R)=0.$ Then the expansion constant of the cochain complex with integer coefficients agrees with the expansion constant with real coefficients.
\end{thm}

Let $G$ be a group of type $\FP_2$ and let $X$ be a 2-complex with $H_1(X;\Z)=0$ such that $G$ acts freely, cellularly, and cocompactly on $X$. Let $F_*\to\Z$ be the associated cellular partial resolution.
For $H\normalsub G$, there is an identification for $i\leq2$ $$C_i(F_*|_H;\Z)\cong C_i(X/H;\Z);$$ and likewise for cochains  $$C^i(F_*|_H;\Z)\cong C^i(X/H;\Z).$$ 
Now, Theorem \ref{linearprogram} applies to $X/H$ when $H$ has (T), as this implies the complex has trivial first cohomology. This now implies the following, exactly as in \cite{BaderSauer} Theorem 2.13.

\begin{prop}\label{preintegralexpansion} Let $X$ be a 2-complex with $H_1(X;\Z)=0$ such that a group $G$ with property (T) acts freely, cellularly, and cocompactly on $X$. Then there is an expansion constant $C$ that applies to all finite index normal subgroups $H$ of $G$. That is, the cochain complexes $C^*(X/H;\Z)$ endowed with the $\ell^1$ norm satisfy uniform linear bounds on the norm of integral primitives of integral coboundaries in degree 2.
\end{prop}

We will now show that this in fact applies to all resolutions by showing the $\Z$-expansion property does not depend on the initial resolution.

\begin{lem}\label{lem:homotopylemma}Let $G$ be a group of type $\FP_2$; let $F_*\to \Z$ be a based partial resolution of type $\FL_2$ and let $P_*\to\Z$ be a partial $\Z G$-resolutions of type $\FP_2$ that is free with a basis in degrees 1 and 2. Then if $F_*\to\Z$ satisfies the $\Z$-coefficient version of Theorem \ref{FP2BaderSauerRexpansion}, then so does $P_*\to\Z.$
\end{lem}
\begin{proof}
The argument is essentially identical to Theorem 3.5 in \cite{HanlonPedroza} and seems to go back to Gersten \cite{Gersten}. Extend the partial resolutions to full projective resolutions.
    Any two projective resolutions are chain homotopy equivalent, so by dualizing, there are cochain maps $f^*:C^*(F_*;\Z) \to C^*(P_*;\Z)$ and $g^*:C^*(P_*;\Z) \to C^*(F_*;\Z)$ such that the compositions $f^*\circ g^*$ and $g^*\circ f^*$ are  cochain homotopic to the identity maps. Let $h^*$ be such a  cochain homotopy, so that 
    $$d h^{i}(\eta) + h^{i+1} (d\eta) = f^i\circ g^i(\eta) - \eta.$$
    
    For $i\in\{1,2\}$, the maps $f^i,g^i,h^i$ are all represented by finite matrices with entries in $\Z G$, depending on the choice of bases. These maps are all bounded in the operator norm; see Lemma 2.7 in \cite{HanlonPedroza}.

    Let $H\normalsub G$ be a finite index normal subgroup. Consider the restricted free resolutions $F_*|_H$ and $P_*|_H$ and corresponding cochain complexes $C^*(F_*|_H;\Z)$ and  $C^*(P_*|_H;\Z)$. The maps $f^*,g^*,h^*$ give cochain maps and a cochain homotopy for the restricted cochain complex and moreover in degrees $i\in\{1,2\},$ these maps have operator norm bounded by a constant $K$ independent of $H$. This is because the induced maps decompose over cosets $Hg\in H\backslash G$ and the corresponding matrices are obtained from the original matrices by replacing the $g\in G$ terms in the entries by permutation matrices and elements of $H$. The calculation in Lemma 2.7 of \cite{HanlonPedroza} now implies the uniform bound.

    In what follows, we drop the superscript notation to unclutter the notation.
    Let $\eta_P$ be an arbitrary coboundary in $C^2(P_*|_H;\Z)$ and set $$\eta_F = g(\eta_P)\in C^2(F_*|_H;\Z).$$ Because $\eta_P$ is exact and $g$ is a cochain map, $\eta_F$ is exact. There is therefore a cochain $\omega_F$ such that $d\omega_F = \eta_F$. Since $F_*$ is assumed to satisfy Theorem \ref{FP2BaderSauerRexpansion} with $\Z$-coefficients, we can assume furthermore that this primitive satisfies $$||\omega_F||\leq C||\eta_F||\leq CK||\eta_P||. $$
    The cochain homotopy condition says $$d h( \eta_P) + h (d\eta_P) = f\circ g(\eta_P) - \eta_P.$$ 
    First observe that $\eta_P$ is coclosed, so $h(d\eta_P)$ vanishes.
     By rearranging, we find $$\eta_P = f\circ g(\eta_P) - dh(\eta_P) = f(\eta_F) - dh(\eta_P).$$
    Set  $\omega_P = f(\omega_F) -h(\eta_P)$ and notice that $d\omega_P = \eta_P.$
Combined with the estimates above, we have $$||\omega_P||\leq K||\omega_F|| + K||\eta_P||\leq CK^2||\eta_P|| + K||\eta_P||,$$ so that indeed $C^*(P_*|_H;\Z)$ satisfies the Theorem \ref{FP2BaderSauerRexpansion} with $\Z$ coefficients and constant $(C
K+1)K$. \end{proof}

\begin{thm}[$\FP_2$ version of \cite{BaderSauer} Theorem 2.13]\label{integralexpansion} Let $G$ be a group of type $\FP_2$ with property (T). Let $F_*\to \Z$ be a based partial $\Z G$-resolution of type $\FL_2$. There is a constant $C$ depending on the partial resolution and bases such that for any finite index normal subgroup $H\normalsub G$ and every coboundary $\eta\in C^2(F_*|_H,\Z)$, there is a cochain $\omega\in C^1(F_*|_H,\Z)$ such that $$d\omega = \eta \text{ and } \norm{\omega}\leq C \norm{\eta}.$$
\end{thm}

\begin{proof}
    This follows from Lemma \ref{overflow}, Proposition \ref{preintegralexpansion} and Lemma \ref{lem:homotopylemma}.
\end{proof}


\subsection{Poincaré duality}
A group $G$ is a Poincaré duality group of dimension $n$, or a $\PD_n$ group, for short, if the following two conditions hold:
\begin{enumerate}
    \item The group $G$ has a finite length resolution by finitely generated projective modules.
    \item  The cohomology $H^i(G;\Z G)$ is concentrated in a single degree $i=n$, and in this degree $H^n(G;\Z G)$ is isomorphic to the possibly nontrivial $\Z G$-module $\Z$.
\end{enumerate}

If $H^n(G;\Z G)$ is isomorphic to $\Z$ with nontrivial module structure, then $G$ is said to be a nonorientable Poincaré duality group. If instead $H^n(G;\Z G)$ is isomorphic to $\Z$ with the trivial $\Z G$-module structure, then we say $G$ is an orientable $\PD_n$ group and write $\PD_n^+$ for short. In the nonorientable case, there is an index two subgroup that is a $\PD_n^+$ group. In the rest of this section, we assume that $G$ is a $\PD_n^+$ group.

Let $F_*\to \Z$ be a partial $\Z G$-resolution of type $\FL_{n-1}$ and then let $0\to F_{n}\to F_*\to \Z$ be a length $n$ projective resolution with $F_{n}$ a finitely generated projective module. Such a projective resolution for $G$ exists by \cite{Brown} Section VIII.6 (this is essentially exercise 2).

Set $F_{i}^\vee:=\hom_{\Z G}(F_{i},\Z G)$; this is naturally a right $\Z G$-module but we twist the group action via the involution $g\mapsto g^{-1}$ to make it a left $\Z G$-module. Because we have a projective resolution of length $n,$ there is a surjection $F_n^\vee\to H^n(G;\Z G).$ 
We can therefore consider the sequence $$F_{0}^\vee\to F_{1}^\vee\to\cdots\to F_{n-1}^\vee\to F_n^\vee\to H^n(G;\Z G).$$ 
The condition that $G$ is a $\PD_n^+$ group ensures the above sequence is a partial projective resolution of the trivial $\Z G$-module $\Z\cong H^n(G;\Z G)$, and the construction ensures it is finitely generated in all degrees and each term $F_{i}^\vee$ is free, except for $F_n^\vee.$ We write $F_{n-*}^\vee\to\Z$ for this new partial resolution.
There are canonical isomorphisms $$\hom_{\Z G}(F_i^\vee,V)\cong \hom_{\Z G}(F_i^\vee,\Z G)\otimes_{\Z G} V \cong F_i\otimes_{\Z G} V$$ induced by the isomorphism $\hom_{\Z G}(F^\vee_i,\Z G)\cong F_i$ and the fact each $F_i$ is projective and finitely generated. Indeed, there is a commutative diagram
\[
\begin{tikzcd}
C^{n-i}(F_{n-*}^\vee;V) \arrow[r,"d"] \arrow[d,leftrightarrow,"\cong"'] &
C^{n-i+1}(F_{n-*}^\vee;V) \arrow[d,leftrightarrow,"\cong"] \\
C_i(F_*;V) \arrow[r,"\d"] &
C_{i-1}(F_*;V)
\end{tikzcd}
\]

\begin{lem}\label{lem:PDisometry} Let $G$ be a $\PD_n^+$ group. Let $F_*\to\Z$ be a projective resolution of $\Z G$-modules of type $\FP_{n}$ extending a based partial resolution of type $\FL_{n-1}.$
 Let $V$ be a normed $\Z G$-module such that $G$ acts linearly and isometrically on $V$. For each $i\neq n,$ give $\hom_{\Z G}(F_i^\vee,V)$ the dual basis defined by the canonical isomorphism and the basis for $F_i$. Then
 the isomorphism $C_i(F_*;V) \cong C^{n-i}(F_{n-*}^\vee;V)$ is an isometry of the normed chain and cochain complexes in degrees $0<i<n$, using the $\ell^1$-norm induced by these bases, as described in the previous section.    
\end{lem}

\begin{proof}
Write $\phi:C_i(F_*;V)\to C^{n-i}(F_{n-*}^\vee;V)$ for the isomorphism described in the commutative diagram above. For $0<i<n$, both the modules $F_i$ and $F_{n-i}$ are free, so the dual basis determines the isomorphism $C^{n-i}(F_{n-*}^\vee;V)\cong V^{m_i}$, which in turn determines the $\ell^1$-norm. The basis also determines an isomorphism $C_i(F_{*};V)\cong V^{m_i}$. Using these identifications, the map  $V^{m_i}\to V^{m_i}$ induced by $\phi$ and these isomorphisms is the identity, thus we have an isometry of $\ell^1$-norms. By the diagram above, this gives an isometry of the normed complexes in the corresponding degrees.
\end{proof}

The above discussion implies the coboundary expansion results from the previous section can be turned into homological expansion in codimension 2.

\begin{prop}\label{pdnhomologicalexpansion}
    Let $G$ be a $\PD_n^+$ group with property (T), where $n>2$. Let $F_*\to\Z$ be a
 $\Z G$-resolution of type $\FP_{n}$ extending a based free resolution of type $\FL_{n-1}$. Then there is a constant $C$ depending on the partial resolution and bases, such that for any finite index normal subgroup $H\normalsub G$ and every exact $(n-2)$-cycle $z\in C_{n-2}(F_*|_H,\Z)$, there is a  $(n-1)$-chain $A\in C_{n-1}(F_*|_H,\Z)$ such that $$\d A = z \text{ and } \norm{A}\leq C \norm{z}.$$
\end{prop}
\begin{proof}
    Let $F_{n-*}^\vee\to\Z$ be the Poincaré dual based partial resolution described above. Let $F'_*\to\Z$ be any based partial resolution of type $\FL_{n-1}$. Apply Theorem \ref{integralexpansion} to $F'_*\to\Z$ to control primitives of coboundaries in $C^{2}(F_*'|_H;\Z)$. Use Lemma \ref{lem:homotopylemma} to transfer the coboundary expansion from $C^{2}(F_*'|_H;\Z)$ to $C^{2}(F_{n-*}^\vee|_H;\Z)$. Then apply the Poincaré duality isomorphism between the complexes $$C_{n-i}(F_*|_H;\Z)\cong C^{i}(F_{n-*}^\vee|_H;\Z)$$ as in Lemma \ref{lem:PDisometry}, to obtain the proposition.
\end{proof}

\section{Hyperbolicity and expansion}\label{sec:hypandexpand}

\subsection{Expansion to hyperbolicity}

In this section, we show that a homological version of the expansion from the previous section implies hyperbolicity. Note that later in Section \ref{nonexpansion}, we show that hyperbolicity obstructs homological expansion in degree two. In the manifold setting, Kielak and Nowak linked coboundary expansion (using the Hamming norm) to hyperbolicity \cite{KielakNowak}. 

First we fix some notation. 
Let $G$ be a group of type $\FP_2$ and let $X$ be a 2-complex with $H_1(X;\Z)=0$ on which $G$ acts freely, cellularly, and cocompactly. Let $F_*\to \Z$ be the associated based partial resolution associated to the augmented chain complex. Denote by $\d:F_i\to F_{i-1}$ the maps in the partial resolution. 

We will measure the complexity of boundaries using filling functions associated to the norm $||\cdot||$ induced by the basis of cells in $F_i $ for $i\leq2$. The filling function measures the minimal norm of a $2$-chain with boundary $z$:
$$
\Fill_{F_1}(z) := \inf\big\{|| A||~:~ \d A= z,~ A \in C_{2}(X;\Z)\big\}.
$$
We will require the following homological characterization of hyperbolicity, due to Gersten.

\begin{thm}[Theorem 5.2 \cite{Gersten}]\label{thm:hyp} If there is a constant $C$ such that for any $z\in \d_2(F_2),$ the filling norm function $\Fill_{F_1}(z)\leq C||z||$, then $G$ is hyperbolic.
\end{thm}

\begin{lem}\label{lem:rf} Let $G$ be a residually finite finitely generated group and let $X$ be a cell complex with a free cocompact cellular $G$-action. Let $Y\subset X$ be a finite subcomplex of $X$. Then for any residual chain $G_i<G$, there is a subgroup $H=G_i$ in the chain such that $Y$ projects injectively to $X/H$.
\end{lem}
\begin{proof} Let $\pi:X\to X/H$ be the quotient map.
    Let $p,q\in Y$. Then $\pi(p)=\pi(q)$ if $p=hq$ for some $h\in H.$ Let $S = \{g\in G~:~gY\cap Y\neq\emptyset\}.$ This is finite because $G$ acts properly discontinuously on $X$. Therefore, by residual finiteness, for sufficiently large $i$, we can take $H=G_i$ such that no nontrivial element of $S$ is in $H$. By construction, the projection map $X\to X/H$ restricted to $Y$ is injective, as desired.
\end{proof}

Let $M$ be a finite cell complex. Define $$\rho(M;\Z) = \inf\limits_{z\in \d C_2(M;\Z)\setminus\{0\}}\frac{||z||}{\Fill_M(z;\Z)},$$ where the filling function is defined as $$\Fill_M(z;\Z) := \inf\big\{|| A||~:~ \d A= z,~ A \in C_{2}(M;\Z)\big\}.$$

Observe that this just encodes the homological version of the (reciprocal of the) expansion constants considered earlier.
By convention, we set $\rho(M;\Z)=0$ if $H_1(M;\Q)$ is nontrivial.
Next we show that duality and coboundary expansion imply hyperbolicity.

\begin{prop}\label{hypfromexpansion}
    Let $G$ be a residually finite group of type $\FP_2$. Let $X$ be a cell complex such that $H_1(X;\Z)=0$ and $G$ acts freely, cellularly, and cocompactly on $X$. Moreover, assume there is a residual chain $G_i\normalsub G$, such that $\inf_i\rho(X/G_i;\Z)>\e$. Then $G$ is hyperbolic. 
\end{prop}

\begin{proof}
    Let $z\in \d C_2(X;\Z)$ be nontrivial and assume $z$ has optimal filling $A\in  C_2(X;\Z)$. Let $Y\subset X$ be the subcomplex consisting of all 2-cells in $X$ that are connected to the support of $z$ by a sequence of at most $R$ 2-cells $\{A_i\}$ with $A_i\cap A_{i+1}\neq\emptyset$, where $R>||A||+1$; note that this contains the support of $A.$
    
   For $i$ sufficiently large, by Lemma \ref{lem:rf} we can take $H=G_i$ such that $Y$ projects injectively to $X/H$. The chain $\pi(A)$ therefore bounds $\pi(z)$ in $C_*(X/H;\Z)$, and $\norm{\pi(A)}=\norm{A}$ and $\norm{\pi(z)}=\norm{z}$ due to the injectivity of the projection map on the set $Y$ containing the support of these chains. 
     
         For any 2-chain $A'$ with boundary $\pi(z)$, the construction of $\pi(Y)$ ensures that either $A'$ has support contained in $\pi(Y)$, or else has norm greater than $||A||.$ 
         
         To see this, note that the chain $A'$ can be decomposed as $A'_0+A'_1$ where $A_0'$ and $A'_1$ have disjoint supports and such that:
         \begin{itemize}
             \item  $A'_0$ has support that is \emph{not} connected to the support of $\pi(z)$ by a sequence of adjacent 2-cells in the support of $A'$,
             \item every cell in the support of $A_1'$ \emph{is} connected to the support of $\pi(z)$ by a sequence of adjacent 2-cells in the support of $A'$.
         \end{itemize}
        Because $A_1'$ and $A_0'$ have disjoint support, $||A'_1||\leq ||A'||.$ As we are interested in optimal fillings, we can throw away $A'_0$ and assume $A'=A'_1.$
        In this case, either $A'$ has support contained in $\pi(Y)$, as claimed, or there is a sequence of adjacent 2-cells running from the support of $\pi(z)$ out of $\pi(Y).$ By construction of $\pi(Y)$, this requires at least $R$-many distinct 2-cells to be in the support of $A_1'$. Since $R>||A||$, we conclude $||A||\leq ||A'||$ as claimed. 
         It follows that $\pi(A)$ is the optimal filling of $\pi(z)$ in $X/H.$
         By hypothesis, the optimal filling of $\pi(z)$ in $X/H$ has norm satisfying $||\pi(A)||\leq\frac{1}{\e}||z||$. Thus $A$ gives a linearly bounded filling of $z$ in $X$, so Theorem \ref{thm:hyp} implies $G$ is hyperbolic.
\end{proof}



\subsection{Hyperbolicity to nonexpansion}\label{nonexpansion}
In this section, we show that residually finite hyperbolic groups cannot be homological expanders. The basic idea in this section is that the shortest homotopically essential loop is contained in a large tube, which forces any filling to be very large.

We work in this section with polyhedral cell complexes. 
Make the 1-skeleton of such a 2-complex a metric space by assigning length one to every edge and then taking the path metric.
For a positive integer $R$, the cellular radius $R$ neighborhood of a subset $Y$ of the 1-skeleton is the full subcomplex spanned by all vertices distance at most $R$ from $Y$ in the 1-skeleton's path metric. In particular, a 2-cell $\sigma$ is contained in the radius $R$ neighborhood of a point $x$ if and only if every vertex in the boundary of $\sigma$ is within radius $R$ of $x$ in the path metric on the 1-skeleton. 
Let $x$ be a vertex and $\cB$ the radius $R$ cellular neighborhood of $x.$ Define the boundary of $\cB$, denoted $\d \cB$, to be full subcomplex spanned by vertices that are in $\cB$ but are also contained in cells not in $\cB.$

\begin{lem}\label{lem:tube}
Let $M$ be a finite 2-complex with fundamental group $G$. Let $g$ be the shortest homotopically essential cellular loop in $M$; let $L$ be the length of $g.$  Fix a lift $\tilde g$ of $g$ to the universal cover $\widetilde M$ and a vertex $x\in \tilde g$. Let $\cB$ be the cellular neighborhood of radius $R\leq \floor{(L-1)/4}$ of $x$ in $\widetilde M$. Then for any nontrivial $h\in G$, $h \cB \cap \cB=\emptyset.$     
\end{lem}
\begin{proof} 

Suppose not, then there are vertices $y,z\in \cB$ such that $hy=z$ for nontrivial $h\in G$. Because $y$ and $z$ are in $\cB$, there are cellular arcs of length at most $R$ connecting $y$ and $z$ to $x$, thus they are connected by a cellular path of length at most $2R$. This path projects to a homotopically nontrivial closed loop that is shorter than $g$, giving a contradiction. \end{proof}

\begin{lem}\label{lem:component}
    Let $M$ be a finite 2-complex. Let $g$ be the shortest homotopically essential cellular loop in $M$. Let $R$ be a positive integer and let $\cB$ be the radius $R$ cellular neighborhood of a vertex $x$ in $g$. Then $g\cap \cB$ has a single component. If $R\leq\floor{(L-1)/4}, $ where $L$ is the length of $g,$ then $g\cap\cB$ is an arc with distinct endpoints.
\end{lem}
\begin{proof}
    First note $g$ is embedded, as otherwise one could construct a shorter homotopically essential loop.
    Suppose that $g\cap \cB$ has multiple components. Let $g_x^\cB$ be the length $2R$ subsegment of $g$ containing $x$ in $\cB$ and let $g_y^\cB$ be another component of $g\cap \cB$ containing a vertex $y$ with $d(x,y)\leq R$. Let $\tau$ be a minimal length path from $x$ to $y$. Let $\alpha,\beta$ be the two subarcs of $g$ connecting $x$ and $y$ extending subsegments of $g_x^\cB$ and oriented to agree with the orientation of $g$. Both of these arcs must leave the radius $R$ cellular neighborhood, so have length at least $R+1$, as they contain distinct vertices distance $R+1$ from $x$. 

    Consider the loops $g_\alpha = \alpha*\tau$ and $g_\beta = \tau^{-1}*\beta$, where $*$ denotes concatenation. Since $g_\alpha*g_\beta$ is freely homotopic to $g$, at least one of $g_\alpha$ and $g_\beta$ is homotopically essential. By construction, both $g_\alpha$ and $g_\beta$ are shorter than the $g$, as we replaced an arc of length at least $R+1$ with an arc whose length is at most $R$. But this contradicts that the initial loop was the shortest homotopically essential loop in $M.$   

    If $R\leq \floor{(L-1)/4}$, then by Lemma \ref{lem:tube}, the ball $\cB$ can be lifted to $\widetilde M$ so that the projection map is injective on the lift. If $g\cap\cB$ is a closed loop, then $g$ is contained in $\cB$ and can therefore be lifted as a closed loop to $\widetilde M$, contradicting that $g$ is homotopically essential. Thus $g\cap\cB$ is an arc with distinct endpoints.
    \end{proof}

\begin{lem}\label{lem:interface}
    Let $M$ be a finite 2-complex. Let $g$ be a cellular loop in $M$ and let $x$ be a vertex in $g$; let $\cB$ be the radius $R$ cellular neighborhood of $x$ for some integer $R$. Let $A$ be 2-chain with boundary $\d A = mg$ and set $A|_\cB$ to be the restriction of $A$ to $\cB.$ Then $\d (A|_\cB)=mg|_\cB + c$,  where $c$ is supported in $\d \cB$.
\end{lem}
\begin{proof}
    Set $c = \d (A|_\cB) - (\d A)|_\cB = \d (A|_\cB) -mg|_\cB.$ Every edge in $\d (A|_\cB)$ is in $\cB,$ because it is in the boundary of a 2-cell that is in $\cB.$  Let $A' = A - A|_\cB.$ Every 2-cell in the support of $A'$ is not in $\cB.$ We have $\d A' = \d A - \d (A|_\cB)=mg - mg|_\cB - c.$ Since $m(g-g|_\cB)$ is supported on edges not in $\cB$ and $c$ is supported on edges that are in $\cB,$ it follows that $c = -(\d A')|_\cB$. Thus $c$ is supported on edges in $\d\cB.$ 
    \end{proof}

\begin{lem}\label{lem:boundary_distance}Let $M$ be a finite 2-complex with universal cover $\widetilde M$. Let $x\in \widetilde M$ be a vertex and let $\cB$ be the radius $R$ cellular neighborhood of $x$. Then there is a constant $R_0$ depending only on $\widetilde M$ such that every vertex in $\d \cB$ is distance at least $R-R_0$ from $x$.
\end{lem}
\begin{proof}
A vertex is in $\d \cB$ if it has distance at most $R$ and is contained in a cell that is not contained in $\cB.$ If the corresponding cell not in $\cB$ is an edge, then the vertex is exactly distance $R$ from $x.$
Suppose instead the corresponding cell is a 2-cell $\sigma$.
Some vertex $v$ of $\sigma$ is distance greater than $R$ from $x.$ Since every vertex in $\sigma$ is within distance $R_0$ of $v$, where $R_0$ is the maximal length of the boundary of a 2-cell in $\widetilde M$ (this is finite as $\widetilde M$ has a finite quotient), it follows that every vertex is at least distance $R-R_0$ to $x.$
\end{proof}

Lastly, we require a basic lemma about hyperbolic spaces.

\begin{lem}[Bridson-Haefliger III.H.1.6]\label{lem:Div} Let $X$ be a $\delta$-hyperbolic geodesic space. Let $c$ be a continuous
rectifiable path in $X$. If $[p,q]$ is a geodesic segment connecting the endpoints of $c$,
then for every $x\in [p,q]$ $$d(x,c)\leq \delta|\log_2(\len(c))| + 1.$$
\end{lem}

We now turn to the geometric estimate that ensures all fillings are large.

\begin{lem}\label{lem:superlinear}
Let $M$ be a finite 2-complex with $\delta$-hyperbolic universal cover and fundamental group $G$. Let $g$ be the shortest homotopically essential cellular loop in $M$ with length $L$. Let $x \in g$ and let $\cB$ be the radius $R$ cellular neighborhood of $x$, where $R=\floor{(L-1)/4}$. Let $g|_\cB = g\cap \cB.$ 
      Let $A$ be an integral 2-chain in $\cB$ such that $\d A = mg|_\cB + c$ for an integer $m$, where $c$ is supported in $\d \cB$. Then there are constants $C$ and $R_0$ depending only on $\widetilde M$ such that $$|m|2^{(R-R_0-1)/\delta}\leq C ||A||+ 2|m|R.$$
\end{lem}

\begin{proof}
    By Lemma \ref{lem:component}, the intersection $g|_\cB$ consists of a single arc with distinct endpoints.
    Let $u$ and $v$ be the starting and endpoints of $g|_\cB$. Observe that $\d c = -\d mg|_\cB = m(u-v).$
    Because $c$ is an integral 1-chain, after choosing a gluing for incident edges in $c$, it is a union of at least $|m|$ paths $\{c_k\}$ joining $u $ and $v$, possibly along with additional loops that we can ignore. Furthermore, note that the gluing construction implies $\sum \len(c_k)\leq  ||c||$.
    
    Lemma \ref{lem:tube} implies the cellular neighborhood $\cB$ can be lifted to $\widetilde M$ so that the projection map restricted to the lift is injective; if $\widetilde x$ is the lift of $x$, then the lifted cellular neighborhood $\widetilde \cB$ is the radius $R$ cellular neighborhood of $\widetilde x$. By Lemma \ref{lem:boundary_distance}, every vertex on each lifted path $\widetilde c_k$ lies at least distance $R-R_0$ from $\widetilde x$. Thus every point on $\widetilde c_k$ does as well.
        
    Observe that the lift $\widetilde{g|_\cB}$ is a geodesic path between its endpoints, as otherwise we could modify a lift $\widetilde g$ of $g$ by replacing the subsegment $\widetilde{ g|_\cB}$ with a shorter arc with the same endpoints, and then project to obtain a shorter loop in $M$ homotopic to $g.$
    
    Therefore, we can apply Lemma \ref{lem:Div} in the 1-skeleton of $\widetilde M$ to estimate $$2^{(R-R_0-1)/\delta}\leq \len(\widetilde c_{k})$$ for each lifted path $\widetilde c_{k}$, then project to conclude $2^{(R-R_0-1)/\delta}\leq \len( c_{k}).$
    There is a uniform bound on the number of 1-cells in the boundary of a 2-cell, which depends only on $\widetilde M.$ This implies there is a constant $C$ such that the $\ell^1$-norm of $A$ satisfies $$\sum_k \len(c_{k})\leq ||c||\leq ||\d A|| + |m|||g|_\cB||\leq  C||A|| + 2|m|R.$$
    
    Because there are at least $|m|$ paths $c_{k}$, it follows that $$|m|2^{(R-R_0-1)/\delta}\leq C ||A|| + 2|m|R.\qedhere$$\end{proof}
    
\begin{lem}\label{lem:nogap}
    Let $M$ be a finite 2-complex with $\delta$-hyperbolic universal cover $\widetilde M$ and fundamental group $G$. For any residual chain $G_i\normalsub G$ corresponding to finite covers $M_i=\widetilde M/G_i\to M$ such that the shortest homotopically essential cellular loop in $M_i$ is rationally nullhomologous, the filling constant $\rho(M_i;\Z)$  tends to zero as $i\to\infty$.
\end{lem}

\begin{proof}
    Let $g_i$ be the shortest homotopically essential cellular loop in $M_i$ and denote the length of $g_i$ by $L_i$. As the $M_i$ form a residual tower of covers, $L_i\to\infty$. Set $R_i= \floor{(L_i-1)/4}$.
    By hypothesis, $g_i$ is rationally nullhomologous, so there exists some positive integer $d_i$ such that $d_ig_i =\d A_i$ where $A_i$ is a 2-chain in $C_2(M_i;\Z)$ with norm $$||A_i||=\Fill_{M_i}(d_ig_i;\Z).$$ 
    Let $\cB_i$ be the radius $R_i$ cellular neighborhood of a vertex $x\in g_i$, as in Lemma \ref{lem:tube}. Let $A_i'$ be the restriction of $A_i$ to the neighborhood $\cB_i$ and let $g_i|_{\cB_i}$ be the restriction of $g_i$ to $\cB_i.$ Then by Lemma \ref{lem:interface}, $\d A_i' = d_ig_i|_{\cB_i}+c$, where $c$ is supported on $\d \cB_i$.
    By Lemma \ref{lem:superlinear} applied to $A_i'$,  we have $$d_i2^{(R_i-R_0-1)/\delta}\leq C||A_i'||+2d_iR_i\leq C ||A_i||, $$ where $R_0$ and $C$  depend only on $\widetilde M$ and the final inequality follows from the fact there are at least $2d_iR$ edges in the support of $d_ig_i$ that are not contained in the ball $\cB_i$, so not counted by the previous estimate. Since $R_i$ grows linearly with $L_i$, dividing the left-hand-side by $d_iL_i=||d_ig_i||$, for $L_i\to\infty
    $ we have $$2^{(R_i-R_0-1)/\delta}/L_i\to\infty$$ and thus  $$\norm{A_i}/(d_iL_i)\to\infty .$$ It follows that $\rho(M_i;\Z)\to0.$
\end{proof}

Proposition \ref{hypfromexpansion} combined with Lemma \ref{lem:nogap} together imply no infinite residually finite group ever has covers with uniform boundary expansion. Note that this equivalent to Theorem \ref{introthm:nonexpansion} in the introduction.

\begin{prop}\label{prop:main}
    Let $M$ be a finite 2-complex with infinite fundamental group $G$. If $M_i\to M$ is a residual tower of covers, then $\rho(M_i;\Z)\to0.$ 
\end{prop}
\begin{proof}
    Suppose not, then $\rho(M_i;\Z)$ is uniformly bounded away from zero. By Proposition \ref{hypfromexpansion}, $G$ is hyperbolic. The shortest homotopically essential cellular loop in $M_i$ is rationally nullhomologous, as otherwise $\rho(M_i;\Z)=0$ contradicting the supposition. Thus the conditions for Lemma \ref{lem:nogap} hold. Therefore Lemma \ref{lem:nogap} implies $\rho(M_i;\Z)\to0$, contradicting the supposition. \end{proof}

We can now prove our main result.

\begin{thm} Any residually finite $\PD_3$ group $G$ does not have property (T).
\end{thm}
\begin{proof}
    Suppose not. Consider a cell complex $X$ such that $H_1(X;\Z)=0$ and $G$ acts freely, cellularly, and cocompactly on $X.$    We can assume $G$ is an oriented $\PD_3$ group after possibly passing to an index two subgroup.
 Apply Proposition \ref{pdnhomologicalexpansion} to a partial resolution extending the based cellular resolution $C_*(X;\Z)\to\Z$. Set $M_i = X/G_i$; this gives uniform integral homological expansion: $\rho(M_i;\Z)>\e$. By Proposition \ref{hypfromexpansion}, $G$ is hyperbolic and therefore finitely presented. 
    
    Let $M$ be a presentation 2-complex. Let $G_i\normalsub G$ be a residual chain and apply Proposition \ref{prop:main} to $\widetilde M/G_i$ to see that for $i\to\infty,$ we have $\rho(\widetilde M/G_i;\Z)\to0$. But this contradicts Proposition \ref{pdnhomologicalexpansion} applied to the cellular resolution associated to the augmented chain complex $C_*(\widetilde M;\Z)$. Thus, no such $G$ exists. 
\end{proof}


\section{Fundamental groups of 3-manifolds and (T)}

In this section we show that 3-manifolds with infinite fundamental group never have (T). The proof is essentially just Lemma \ref{lem:nogap} combined with the expansion result of Bader and Sauer and chain-level Poincaré duality applied to triangulations (which is a key tool in the waist inequalities in \cite{BaderSauer}).

\begin{lem}\label{lem:homologicalgap}
    Let $M$ be an oriented closed 3-manifold with fundamental group $G$ with property (T). 
    Fix a triangulation $\cT$ of $M$ with dual cellulation $\cT^\vee$. Let $M_i$ be any sequence of regular finite covers of $M$ with pullback dual cellulations $\cT_i^\vee$. Then $\rho(\cT^\vee_i;\Z)>\e$ for some $\e$ depending only on $\cT.$ 
\end{lem}
\begin{proof}
    There is an $\ell^1$-isometric chain map $\phi: C_{3-*}(\cT_i^\vee;\Z)\to C^*(\cT_i;\Z)$ inducing the Poincaré duality isomorphism on homology. Theorem \ref{integralexpansion} and this chain-level Poincaré duality imply a uniform lower bound on $\rho(\cT^{\vee}_i;\Z).$
\end{proof}

\begin{thm} Let $G$ be the fundamental group of a compact 3-manifold $M$. If $G$ has property (T), then $G$ is finite.
\end{thm}
\begin{proof}
    We first can assume $M$ is orientable, as if it is not, we just replace $M$ with its orientation double cover.
    We can also assume no boundary component is a sphere, as attaching 3-balls along spheres does not change the fundamental group. If any boundary component is not a sphere, then $H^1(M;\Z)$ is nontrivial by the "half lives, half dies" lemma and Poincaré-Lefschetz duality, which is incompatible with property (T). Thus we can assume $M$ is a closed orientable 3-manifold.
    
    The fundamental groups of 3-manifolds are residually finite by an argument of Hempel \cite{Hempel} combined with the geometrization theorem of Perelman \cite{PerelmanEntropy}. Fix a closed 3-manifold $M$ with infinite fundamental group with (T). Take a residual tower of covers $M_i\to M$ and apply Lemma \ref{lem:homologicalgap} to a triangulation of $M$, giving each cover $M_i$ the pullback triangulation. This now contradicts Proposition \ref{prop:main}, applied to the 2-skeleton of the dual cell complex. Thus no such $M$ exists. 
\end{proof}

\subsection*{Acknowledgments}
The author thanks Marc Lackenby for helpful conversations and Shaked Bader for helpful comments on a draft.
\small{
\bibliographystyle{alpha}
\bibliography{bib}}

\end{document}